\newcommand{\al}{\alpha}               \newcommand{\be}{\beta}
\newcommand{\ga}{\gamma}               \newcommand{\Ga}{\Gamma}
\newcommand{\de}{\delta}               \newcommand{\De}{\Delta}
\newcommand{\lb}{\lambda}              \newcommand{\Lb}{\Lambda}
              
\newcommand{\veps}{\varepsilon}        \newcommand{\vphi}{\varphi}

\newcommand{\cal}{\mathcal}
\newcommand{\calb}{{\cal B}}           \newcommand{\calc}{{\cal C}}
             
\newcommand{\calh}{{\cal H}}           \newcommand{\cali}{{\cal I}}
\newcommand{\calp}{{\cal P}}           \newcommand{\calr}{{\cal R}}
\newcommand{\cals}{{\cal S}}           
\newcommand{\calu}{{\cal U}}           \newcommand{\calv}{{\cal V}}

\newcommand{\Dom}{{\rm Dom}}           
\newcommand{\cl}{{\rm cl}}             \newcommand{\gr}{{\rm gr}}

\newcommand{\incl}{\subseteq}          \newcommand{\aincl}{\supseteq}
       
\newcommand{\es}{\emptyset}            \newcommand{\sm}{\setminus}
\newcommand{\impl}{\Rightarrow}        \newcommand{\limpl}{\Longrightarrow}
    \newcommand{\lequi}{\Longleftrightarrow}
\newcommand{\oo}{\infty}               
             \newcommand{\wt}{\widetilde}

\newcommand{\barr}{\begin{array}}        \newcommand{\earr}{\end{array}}
\newcommand{\beq}{\begin{equation}}      \newcommand{\eeq}{\end{equation}}
\newcommand{\bit}{\begin{itemize}}       \newcommand{\eit}{\end{itemize}}
\newcommand{\blemma}{\begin{lemma}}      \newcommand{\elemma}{\end{lemma}}
\newcommand{\bproof}{\begin{proof}}      \newcommand{\eproof}{\end{proof}}
\newcommand{\bprop}{\begin{proposition}} \newcommand{\eprop}{\end{proposition}}
\newcommand{\brem}{\begin{remark}}       \newcommand{\erem}{\end{remark}}
\newcommand{\btab}{\begin{tabular}}      \newcommand{\etab}{\end{tabular}}
\newcommand{\btheorem}{\begin{theorem}}  \newcommand{\etheorem}{\end{theorem}}

               \newcommand{\sk}{\smallskip}
\newcommand{\n}{\noindent}

\def\Roo{R\cup \{\oo\}}  \def\R+oo{R_+\cup\{\oo\}}

\def\dtends    {\stackrel {\it d}{\longrightarrow}}

\def\dlbtends  {\stackrel {\it d_\lb}{\longrightarrow}}
\def\elbtends  {\stackrel {\it e_\lb}{\longrightarrow}}
\def\Btends    {\stackrel {\calb}{\longrightarrow}}
\def\Ctends    {\stackrel {\calc}{\longrightarrow}}
\def\Dtends    {\stackrel {\it D}{\longrightarrow}}
\def\Etends    {\stackrel {\it E}{\longrightarrow}}
\def\Vtends    {\stackrel {\calv}{\longrightarrow}}
\def\Detends   {\stackrel {\De}{\longrightarrow}}

\documentclass[reqno]{amsart}

\newtheorem{theorem}{\bf Theorem}

\newtheorem{lemma}{\bf Lemma}
\newtheorem{proposition}{\bf Proposition}
\newtheorem{remark}{\bf Remark}

\begin{document}

\title
[VARIATIONAL PRINCIPLES IN FANG UNIFORM SPACES]
{VARIATIONAL PRINCIPLES \\
IN FANG UNIFORM SPACES}

\author{Mihai Turinici}
\address{
"A. Myller" Mathematical Seminar;
"A. I. Cuza" University;
700506 Ia\c{s}i, Romania
}
\email{mturi@uaic.ro}

%\date{February 6, 2013}



\begin{abstract}
The vectorial Zhu-Li Variational Principle (ZLVP) in Fang uniform spaces 
is in the logical segment between 
the Brezis-Browder ordering principle (BB) 
and Ekeland's Variational Principle (EVP);
hence, it is equivalent with both BB and EVP.
In particular, the conclusion is applicable to
Hamel's Variational Principle (HVP). 
Finally, a proof of HVP $\lequi$ EVP is provided, 
by means of a direct approach.
\end{abstract}

\maketitle

% Section 1
\section{Introduction}
\setcounter{equation}{0}

Let $X$ be a nonempty set.
Take a metric $d:X\times X\to R_+$ over it,
where $R_+:=[0,\oo[$;
and let $\vphi:X\to R$ be a function taken as
\bit
\item[(a01)]
$\vphi$ is bounded from below:
[$\vphi(x)\ge b$, $\forall x\in X$], for some $b\in R$.
\eit
The following Ekeland's variational principle 
\cite{ekeland-1979}
(in short: EVP) is our starting point. Assume that (in addition)
\bit
\item[(a02)]
$d$ is complete:
each $d$-Cauchy sequence is $d$-convergent
\item[(a03)]
$\vphi$ is $d$-lsc:
$\liminf_n \vphi(x_n)\ge \vphi(x)$, whenever $x_n\dtends x$.
\eit

% Theorem 1
\btheorem \label{t1}
Let these conditions hold; and $u\in X$ be arbitrary fixed.
There exists then $v=v(u)\in X$  in such a way that
\beq \label{101}
\mbox{
$d(u,v)\le \vphi(u)-\vphi(v)$ (hence $\vphi(u)\ge \vphi(v)$)
}
\eeq
\beq \label{102}
\mbox{
$d(v,x)> \vphi(v)-\vphi(x)$, \quad for all\ $x\in X\sm \{v\}$.
}
\eeq
\etheorem
As a matter of fact, the original result is with
$\vphi:X \to \Roo$ being, in addition, 
{\it proper} ($\Dom(\vphi)\ne \es$).
But, the author's conclusion is obtainable from this restricted version;
just apply Theorem \ref{t1} to the triplet 
$(X(u,\le);d;\psi)$, where $u\in \Dom(\vphi)$,
$X(u,\le):=\{x\in X; u\le x\}$,  $\psi:=\vphi-b$, 
and $(\le)$ is the {\it quasi-order} 
(i.e.: reflexive and transitive relation) described as
\bit
\item[(a04)]
($x,y\in X$):\ $x\le y$ iff $d(x,y)+\vphi(y)\le \vphi(x)$.
\eit

This principle found some basic applications to control and optimization,
critical point theory and global  analysis.  
As a consequence, many extensions of it were proposed;
see, for instance,
Hyers, Isac and Rassias \cite[Ch 5]{hyers-isac-rassias-97}.
Here, we are interested in 
the {\it structural} generalizations of EVP, related to  
the metrical structure 
being substituted by a {\it uniform} one.
The basic contribution 
is Hamel's variational principle 
\cite{hamel-2005}
(in short: HVP);
for "inductive" extensions of it, we refer to the  
Zhu-Li variational principle \cite{zhu-li-2007}
(in short: ZLVP).

Now, it is our aim in the following to show 
(cf. Section 4) that ZLVP is nothing but 
an equivalent version of EVP.
The basic tools of our investigations are 
a lot of countable maximal statements 
(given in Section 2) deductible from 
the (Bernays-Tarski) Principle of Dependent Choices
(in short: DC)
and the concept of conical gauge function 
(developed in Section 3).
In particular (cf. Section 5), the 
logical equivalence above remains valid  
when ZLVP is substituted by HVP.
Finally, in Section 6, 
a direct approach for HVP $\lequi$ EVP is proposed,
without passing through DC.
Some other aspects will be discussed elsewhere.

%  Section 2
\section{Countable maximal statements}
\setcounter{equation}{0}

Let $M$ be a nonempty set; and $\calr\incl M\times M$ stand for 
a (nonempty) relation over it.
For each $x\in M$, denote $M(x,\calr)=\{y\in M; x\calr y\}$.
The following "Dependent Choices Principle" (in short: DC) 
is our starting point:

% Proposition 1 
\bprop \label{p1}
Suppose that
\bit
\item[(b01)]
$M(c,\calr)$ is nonempty, for all $c\in M$.
\eit
Then, for each $a\in M$ there exists $(x_n)\incl M$ 
with $x_0=a$ and $x_n\calr x_{n+1}$, for all $n$.
\eprop

This principle, due to 
Bernays \cite{bernays-1942}
and
Tarski \cite{tarski-1948},
is deductible from AC (= the Axiom of Choice), but not conversely; 
cf. 
Wolk \cite{wolk-1983}.
Moreover, 
the {\it alternate Zermelo-Fraenkel system} (ZF-AC+DC)
seems to be sufficient for a large part of the
"usual" mathematics; see
Moore \cite[Appendix 2]{moore-1982}.
\sk

{\bf (A)}
Let $M$ be some nonempty set. Take a
{\it quasi-order} $(\le)$ over it,  
as well as a function $\vphi:M \to R$.
Call the point $z\in M$, $(\le,\vphi)$-{\it maximal} when:
$\vphi$ is constant on $M(z,\le)$.
A basic result about the existence of such points is the
Brezis-Browder ordering principle \cite{brezis-browder-1976}
(in short: BB).

% Proposition 2
\bprop  \label{p2}
Suppose that
\bit
\item[(b02)]
$(M,\le)$ is sequentially inductive:\\
each ascending sequence has an upper bound (modulo $(\le)$)
\item[(b03)]
$\vphi$ is bounded from below and $(\le)$-decreasing\ ($x\le y \limpl \vphi(x) \ge \vphi(y)$).
\eit
Then, for each $u\in M$ there exists a $(\le, \vphi)$-maximal
$v\in M$ with $u\le v$.
\eprop

\bproof
Define the function $\be:M\to R$ as:
$\be(v):=\inf[\vphi(M(v,\le))]$, $v\in M$.
Clearly, $\be$ is increasing and 
[$\vphi(v)\ge \be(v)$, $\forall v\in M$].
Further, (b03) gives at once a characterization like:
$v$ is $(\le,\vphi)$-maximal iff $\vphi(v)=\be(v)$.
Now, assume by contradiction that 
the conclusion in this statement is false; 
i.e. (see above)
there must be some $u\in M$ such that:
\bit
\item[]
for each $v\in M_u:=M(u,\le)$, one has $\vphi(v)> \be(v)$.
\eit
Consequently (for all such $v$), 
$\vphi(v)> (1/2)(\vphi(v)+\be(v)) > \be(v)$;
hence 
\beq \label{201}
\mbox{
$v\le w$ and $(1/2)(\vphi(v)+\be(v))> \vphi(w)$, 
}
\eeq
for at least one $w$ (belonging to $M_u$).
The relation $\calr$ over $M_u$ introduced {\it via} (\ref{201})
fulfills 
$M_u(v,\calr)\ne \es$, for all $v\in M_u$.
So, by (DC), there must be 
a sequence $(u_n)$ in $M_u$ with $u_0=u$ and
\beq \label{202}
\mbox{
$u_n\le u_{n+1}$, $(1/2)(\vphi(u_n)+\be(u_n))> \vphi(u_{n+1})$, 
for all $n$.
}
\eeq
We have thus constructed an ascending sequence  $(u_n)$ in $M_u$
for which $(\vphi(u_n))$ is
(strictly) descending and bounded below; hence 
$\lb:=\lim_n \vphi(u_n)$ exists in $R$.
Moreover, from (b02), $(u_n)$ is bounded above in $M$: 
there exists  $v\in M$ such that $u_n\le v$, $\forall n$.
Combining with (b03), gives $\vphi(u_n)\ge \vphi(v)$, $\forall n$; 
and (by the properties of $\be$)
$\vphi(v)\ge \be(v)\ge \be(u_n)$, $\forall n$.
The former of these relations gives $\lb\ge \vphi(v)$ 
(passing to limit as $n\to \oo$).
On the other hand, the latter of these relations yields
({\it via} (\ref{202}))
$(1/2)(\vphi(u_n)+\be(v))> \vphi(u_{n+1})$, for all $n\in N$. 
Passing to limit as $n\to \oo$  yields
$(\vphi(v)\ge )\be(v)\ge \lb$; so, combining with the preceding relation,
$\vphi(v)=\be(v)(=\lb)$, contradiction.
Hence, our working assumption cannot be accepted, 
and the conclusion follows.
\eproof

{\bf (B)}
This principle, including (EVP) (see below) 
found some useful applications to convex and nonconvex analysis.
For this reason, it was the subject
of  many extensions; see, e.g., 
Kang and Park \cite{kang-park-1990}.
However, we must emphasize that, whenever a
maximal principle of this type is to be applied, a substitution of it
by the Brezis-Browder's (BB) is always possible.
This raises the question of to what extent are these enlargements of
BB effective. 
Before giving a complete (negative) answer to this, we note that
a way of obtaining structural extensions from BB
is by "splitting" the key condition (b02) as
\bit
\item[(b04)]
($\forall (x_n)\incl M$)
ascending $\impl$ Cauchy $\impl$ 
convergent $\impl$ bounded above.
\eit
This will necessitate some conventions and auxiliary facts.
Let $\cals(M,\le)$ stand for the class of all ascending sequences in $M$.
By a (sequential) {\it convergence structure} on $(M,\le)$ we mean, as in
Kasahara \cite{kasahara-1976},
any part  $\calc$ of $\cals(M,\le)\times M$ with
\bit
\item[(b05)]
$x_n=x, \forall n\in N  \limpl ((x_n);x)\in \calc$
\item[(b06)]
\mbox{
$((x_n);x)\in \calc \limpl ((y_n);x)\in \calc$, for each
subsequence $(y_n)$ of $(x_n)$.
}
\eit
In this case, $((x_n);x)\in \calc$ writes  $x_n \Ctends x$;
and reads: $x$ is the {\it $\calc$-limit} of $(x_n)$.
When such elements exist,
we say that $(x_n)$ is $\calc$-{\it convergent};
the class of all these will be denoted $\cals_\calc(M,\le)$.
Further, by a (sequential) {\it Cauchy structure} on $(M,\le)$ we
mean any part $\calh$  of $\cals(M,\le)$ with
\bit
\item[(b07)]
$x_n=x, \forall n\in N  \limpl (x_n)\in \calh$
\item[(b08)]
$(x_n)\in \calh \limpl (y_n)\in \calh$, for each
subsequence $(y_n)$ of $(x_n)$.
\eit
Each element of $\calh$ will be referred to as a
{\it $\calh$-Cauchy} sequence.
[For example, a good choice is
$\calh=\cals_\calc(M,\le)$; but this is not the only
possible one].
Suppose that we introduced such a couple $(\calc,\calh)$,
referred to as a {\it conv-Cauchy structure}.
Roughly speaking, the objective to be attained is the realization of (b04).
To this end, the following conditions will be considered
\bit
\item[(b09)]
$\calh$ is $(\le)$-regular:
each ascending sequence in $M$ is $\calh$-Cauchy
\item[(b10)]
$(\calc,\calh)$ is (sequentially) $(\le)$-complete:\\
each ascending $\calh$-Cauchy sequence in $M$ is $\calc$-convergent
\item[(b11)]
$(\le)$ is $\calc$-selfclosed:
the $\calc$-limit of each ascending $\calc$-convergent
sequence in $M$ is an upper bound of it.
\eit
The following structural version of BB is then available.
Let again $(M,\le)$ be a quasi-ordered structure;
and $\vphi:M\to R$ be as in (b03).

% Proposition 3
\bprop \label{p3}
Suppose that the conv-Cauchy structure $(\calc,\calh)$ is 
such that (b09)-(b11) hold.
Then, conclusion in BB is retainable.
\eprop

The proof is immediate ({\it via} BB);
just note that (b09)-(b11) imply (b02). Hence, Proposition \ref{p3}
is deductible from BB. The reciprocal deduction
is also possible. To verify this, it is enough to take the
convergence structure over $(M,\le)$ as
the {\it bounded from above} property $\calb$
[introduced as: $x_n \Btends x$\ iff $x_n\le x$, for all $n$];
and the Cauchy structure on $(M.\le)$
be identical with $\cals_\calb(M,\le)$.
\sk

{\bf (C)}
A basic application of these facts is to be done in
the pseudo-uniform setting.
Let $(M,\le)$ be a quasi-ordered structure.
Denote
$\cali(M)=\{(x,x); x\in M\}$ (the {\it diagonal} of $M$);
and let $\calv$ be a family of parts in $M\times M$.
Under a convention similar to that in
Nachbin  \cite[Ch 2, Sect 2]{nachbin-1965},
we say that $\calv$ is a {\it pseudo-uniformity} over it when
$\cap \calv\aincl \cali (M)$.
Suppose that we introduced such an  object.
The associated (sequential) convergence structure
$(\calv)$ on $(M,\le)$ may be described as
\bit
\item[]
$x_n \Vtends x$ iff
$\forall V\in \calv, \exists n(V):\ n\ge n(V)\limpl (x_n,x)\in V$.
\eit
It will be referred to as: $x$ is the $\calv$-limit of $(x_n)$; 
if such elements exist, we say that $(x_n)$ is $\calv$-convergent.
In addition, we may introduce the  $\calv$-{\it Cauchy}
property for an ascending sequence $(x_n)$ as
\bit
\item[]
$\forall V\in \calv,\exists n(V): n(V)\le p\le q\limpl (x_p,x_q)\in V$;
\eit
the class of all these will be denoted as Cauchy($\calv$).
Now, in this context, further interpretations of the
regularity conditions above are possible.
Call the (ascending) sequence $(x_n)$, $\calv$-{\it asymptotic} provided
\bit
\item[]
$\forall V\in \calv, \exists n(V): n\ge n(V)\limpl (x_n,x_{n+1})\in V$.
\eit
Clearly, each $\calv$-Cauchy (ascending) sequence is $\calv$-asymptotic too.
The converse is also true, if {\it all} such sequences are involved;
so that, the global conditions below 
\bit
\item[(b12)]
each ascending sequence in $M$ is $\calv$-Cauchy
\item[(b13)]
each ascending sequence in $M$ is $\calv$-asymptotic
\eit
are equivalent to each other. 
By definition, either of these will be referred to as
$\calv$ is {\it $(\le)$-regular};
this is just (b09), relative to $\calh$=Cauchy($\calv$).
Further, call $\calv$,
{\it (sequentially) $(\le)$-complete} when each ascending
$\calv$-Cauchy sequence in $M$ is $\calv$-convergent.
As before, this is nothing but the condition (b10),
relative to $\calc=(\calv)$ and $\calh$=Cauchy($\calv$).
Finally, let us say that $(\le)$ is
{\it $\calv$-selfclosed} when the $\calv$-limit of each
ascending $\calv$-convergent sequence in $M$
is an upper bound of it; i.e., (b11) holds relative to $\calc=(\calv)$.
The following "uniform" type version of Proposition \ref{p3}
is then available.
[The general conditions about $(M,\le)$ and $\vphi$
prevail].

% Proposition 4
\bprop \label{p4}
Assume that 
$\calv$ is  $(\le)$-regular, (sequentially) $(\le)$-complete
and 
$(\le)$ is $\calv$-selfclosed.
Then, conclusions of BB are retainable.
\eprop

As a consequence of this, 
BB $\limpl$ Proposition \ref{p3} $\limpl$ Proposition \ref{p4}; 
in addition, Proposition \ref{p4} $\limpl$ BB.
In fact, let the premises of BB hold; 
and put $\calv=\{\gr(\le)\}$, where
$\gr(\le):=\{(x,y)\in M\times M; x\le y\}$.
It is easy to see that all conditions in
Proposition \ref{p4} are fulfilled; hence the claim.
Summing up, BB, Proposition \ref{p3} and Proposition \ref{p4} 
are mutually equivalent. 

The discussed particular case is an "extremely"  one.
To get "standard" examples in the area, we need further conventions.
Let $\calv$ be a family of parts in $M\times M$;
we call it a {\it fundamental system of entourages} for a uniformity
over $M$, when
(cf. Bourbaki \cite[Ch 2, Sect 1]{bourbaki-1989})
\bit
\item[(b14)]
$(\calv,\aincl)$ is directed and $\cap\calv \aincl \cali(M)$
\item[(b15)]
$\forall V\in \calv, \exists W\in \calv:\ W\incl V^{-1}, W\circ W\incl V$.
\eit
The uniformity in question is just
$\calu$=\{$P\incl M\times M$; $P\aincl Q$, for some $Q\in \calv$\}.
As a rule, the "uniform" terminology refers to it.
However (as results directly by definition), all $\calu$-notions
are in fact $\calv$-notions.
For example, the (sequential) convergence structure $(\calu)$
on $(M,\le)$
(introduced as before) is nothing else than  $(\calv)$.
Likewise, the (attached to $\calu$) Cauchy and asymptotic properties
are identical with those related to $\calv$.

The following "standard" version of Proposition \ref{p4}
holds. (As before, $(M,\le)$ is a quasi-ordered structure; and
$\vphi:M\to R$ is as in (b03)).

% Proposition 5
\bprop \label{p5}
Assume
$\calv$ is (sequentially) $(\le)$-complete,
$(\le)$ is $\calv$-selfclosed and
\bit
\item[(b16)]
($\calv$ is $(\le,\vphi)$-compatible)
$\forall V\in \calv, \exists \de=\de(V)> 0$:\\
$x,y\in M, x\le y, \vphi(x)-\vphi(y)< \de \limpl (x,y)\in V$.
\eit
Then, for each $u\in M$ there exists a $(\le,\vphi)$-maximal
$v\in M$ with $u\le v$.
\eprop

\bproof
It is sufficient to prove that
$\calv$ is  $(\le)$-regular (see above).
Let $(x_n)$ be an ascending (modulo $(\le)$) sequence
in $M$. The sequence $(\vphi(x_n))$ is descending and bounded
from below; hence a Cauchy one
$$
\forall \de> 0,\ \exists n(\de):\ 
n(\de)\le p\le q \limpl \vphi(x_p)-\vphi(x_q)< \de.
$$
This, along with (b16), gives us the conclusion we want.
\eproof

Note that, a direct consequence of (b16) is
\beq \label{203}
\mbox{
if $x<>y$ and  $\vphi(x)=\vphi(y)$ then $(x,y)\in \cap \calv$.
}
\eeq
(Here, $x<>y$ means: either $x\le y$ or $y\le x$).
This gives 
\beq \label{204}
\mbox{
($\forall z\in M$)
$(\le,\vphi)$-maximal $\limpl$ $(\le,\calv)$-maximal;
}
\eeq
where the last property means:
$w\in M$ and $z\le w$ imply $(z,w)\in \cap \calv$.
So, Proposition \ref{p5} is, at the same time, an existence principle for
$(\le, \calv)$-maximal elements; and, as such, it may be compared
with a related statement of
Turinici \cite{turinici-1994}.
In particular, when $\calv$ is {\it separated} ($\cap \calv =\cali(M)$),
we have (again {\it via} (b16))
\beq \label{205}
\mbox{
($\forall z\in M$):\
$(\le,\calv)$-maximal $\lequi$ $(\le)$-maximal ($M(z,\le)=\{z\})$;
}
\eeq
and Proposition \ref{p5} yields the maximal principle in
Hamel \cite[Theorem 1]{hamel-2005}.
On the other hand,
(\ref{203}) (and the separated property of $\calv$) also gives
({\it via} (b03))
\beq \label{206}
\mbox{
$(\forall x,y\in M)$:\  $x\le y$,\ $y\le x$ $\limpl$ $x=y$;
}
\end{equation}
whence, $(\le)$ is {\it antisymmetric} (hence an {\it order}) on $X$. 
The "separated" variant of Proposition \ref{p5}
is then identical with the main result in
Br{\o}ndsted \cite[Theorem 1]{brondsted-1974};
for this reason, it will be referred to as
the Br{\o}ndsted Maximal Principle (in short: BMP).
Further aspects may be found in
Mizoguchi \cite{mizoguchi-1990}.
\sk

{\bf (D)}
Let VP stand for any of the variational principles 
described by Proposition \ref{p2} -- Proposition \ref{p5}.
Note that, by the developments above, one has
(DC) $\limpl$ (BB) $\limpl$ (VP);
moreover, it is easy to see that (VP) $\limpl$ (EVP).
This raises the question of to what extent 
are these inclusions effective; note that, its natural setting is
(ZF-AC) (the {\it reduced} Zermelo-Fraenkel system).
The answer is negative; 
for a description of it, we need some preliminary facts. 

Let $X$ be a nonempty set; 
and $(\le)$ be an order on it. We say that $(\le)$
has  the {\it inf-lattice} property, provided: 
$x\wedge y:=\inf(x,y)$ exists, for all $x,y\in X$.
Further, we say that $z\in X$ is a $(\le)$-{\it maximal} element if 
$X(z,\le)=\{z\}$; the class of all these points will be
denoted as $\max(X,\le)$. 
In this case, $(\le)$ is called a {\it Zorn order} when 
$\max(X,\le)$ is nonempty and {\it cofinal} in $X$
[for each $u\in X$ there exists a $(\le)$-maximal 
$v\in X$ with $u\le v$].
Further aspects are to be described in a metric setting. 
Let $d:X\times X\to R_+$ be a metric over $X$;
and $\vphi:X\to R_+$ be some function.
Then, the natural choice for $(\le)$ above is
\bit
\item[]
$x\le_{(d,\vphi)} y$ iff $d(x,y)\le \vphi(x)-\vphi(y)$;
\eit
referred to as the
Br{\o}ndsted order \cite{brondsted-1976}
attached to $(d,\vphi)$. 
Denote
$X(x,\rho)=\{u\in X; d(x,u)< \rho\}$, $x\in X$, $\rho> 0$
[the open sphere with center $x$ and radius $\rho$].
Call the ambient metric space  $(X,d)$, {\it discrete} when for
each $x\in X$ there exists $\rho=\rho(x)> 0$ such that 
$X(x,\rho)=\{x\}$. Note that, under such an assumption, 
any function $\psi:X\to R$ is continuous over $X$.
However, the Lipschitz property 
($|\psi(x)-\psi(y)|\le L d(x,y)$, $x,y\in X$, for some $L> 0$)
cannot be assured, in general.

% Theorem 2
\btheorem \label{t2}
Let  the metric space $(X,d)$ and the function $\vphi:X\to R_+$
satisfy
\bit
\item[(b17)]
$(X,d)$ is discrete bounded and complete
\item[(b18)]
$(\le_{(d,\vphi)})$ has the inf-lattice property
\item[(b19)]
$\vphi$ is $d$-nonexpansive and $\vphi(X)$ is countable.
\eit
Then, $(\le_{(d,\vphi)})$ is a Zorn order.
\etheorem

We refer to this statement as: 
the discrete Lipschitz countable version of EVP 
(in short: (EVPdLc)).
Clearly, (EVP) $\limpl$ (EVPdLc). 
The remarkable fact to be added is contained in

% Proposition 6
\bprop \label{p6}
We have (in the reduced Zermelo-Fraenkel system)
(EVPdLc) $\limpl$ (DC).
So (by the above), the maximal/variational principles
(BB), (VP), and (EVP) are all equivalent to (DC);
hence, mutually equivalent.
\eprop

For a detailed proof, see
Turinici \cite{turinici-2011}.
In particular, when the assumptions (b18) and (b19) are
ignored in Theorem \ref{t2}, 
Proposition \ref{p6} reduces to the result in 
Brunner \cite{brunner-1987}.
Further aspects may be found in 
Schechter \cite[Ch 19, Sect 19.53]{schechter-1997}.

% Section 3
\section{Conical gauge functions}
\setcounter{equation}{0}

Let $Y$ be a (real) {\it vector space}.
Take a convex cone $H$ of $Y$
($\al H+\be H\incl H$, for each $\al,\be$ in $R_+$);
which in addition is non-degenerate ($H\ne \{0\}$),
proper ($H\ne Y$);
and let $(\le)$ stand for its induced quasi-order
[$x\le y$ iff $y-x\in H$].
Further, take some point $k^0\in H\sm (-H)$; and put (for $y\in Y$) 
\bit
\item[(c01)]
$\Ga(H;k^0;y)=\{s\in R_+; k^0s\le y\}$,\ $\ga(H;k^0;y)=\sup\Ga(H;k^0;y)$.
\eit
(Here, by convention, $\sup(\es)=-\oo$).
We therefore defined a multivalued function $\Ga(.):=\Ga(H;k^0;.)$
from $Y$ to $\calp(R_+)$, 
and a function $\ga(.):=\ga(H;k^0;.)$ from $Y$ to 
$R\cup\{-\oo\}\cup \{\oo\}$ with 
\beq \label{301}
[\Ga(y)=\es,\ \ga(y)=-\oo] \lequi y\in Y\sm H;
\eeq
the latter of these will be referred to as 
the {\it gauge} function attached to $(H;k^0)$.
Note that, for each $y\in H$,
\beq \label{302}
\mbox{
$\Ga(y)$ is hereditary\ ($s\in \Ga(y)\limpl [0,s]\incl \Ga(y)$);
}
\eeq
so, either 
$\Ga(y)=[0,\ga(y)[$ if $0< \ga(y)\le \oo$, 
or 
$\Ga(y)=[0,\ga(y)]$, if $0\le \ga(y)< \oo$.
In addition, the couple $(\Ga,\ga)$ is 
positively homogeneous and increasing
\beq \label{303}
\Ga(ty)=t \Ga(y), \ga(ty)=t\ga(y),\ \  \forall t> 0,\ \forall y\in H
\eeq
\beq \label{304}
\mbox{
$y_1,y_2\in H$, $y_1\le y_2$\ implies\  
$\Ga(y_1)\incl \Ga(y_2)$, $\ga(y_1)\le \ga(y_2)$.
}
\eeq

An important question to be solved is that of $\Ga$ being
{\it proper} [$\Ga(y)\ne R_+$, $\forall y\in H$].
After 
Cristescu \cite[Ch 5, Sect 1]{cristescu-1977},
we say that $H$ is {\it Archimedean}, provided
\bit
\item[(c02)]
[$v\in Y$, $h\in H$, $\Ga(H;v;h)=R_+$] imply $v\in -H$.
\eit
Likewise, let us say that $H$ is {\it semi-Archimedean}, if
\bit
\item[(c03)]
$\Ga(H;k;y)$ is closed, $\forall k\in H\sm (-H)$,\ $\forall y\in H$.
\eit

% Lemma 1
\blemma \label{le1}
The following are valid:

{\bf i)}
If $H$ is Archimedean, then 
$\Ga(.)$ is proper, in the sense:
$0\le \ga(y)< \oo$ and $\Ga(y)=[0,\ga(y)]$,\ for all $y\in H$;
so, $H$ is semi-Archimedean too

{\bf ii)}
Let $H$ be semi-Archimedean; and
$\al\in R_+$, $y\in H$, $(\be_n; n\ge 0)\incl R_+$ be such that
[$k^0 \al\le y+ k^0 \be_n$, $\forall n$] and  $\be_n\to 0$.
Then, $k^0 \al\le y$.
\elemma

\bproof
{\bf i)}
Let $y\in H$ be arbitrary fixed. 
If $\Ga(y)=R_+$ then, by the Archimedean property of $H$, 
one gets $k^0\in -H$; contradiction.
Consequently, $R_+\sm \Ga(y)\ne \es$;
so that, by (\ref{302}), $\Ga(y)$ is bounded
[whence, $0\le \ga(y)< \oo$].
Further, again by (\ref{302}), 
$k^0\ga(y)-y\le k^0t$, for all $t> 0$; wherefrom
$\Ga(H;k^0\ga(y)-y;k^0)=R_+$.
This, again by the Archimedean property of $H$, gives
$\ga(y)\in \Ga(y)$; i.e., $\Ga(y)=[0,\ga(y)]$.

{\bf ii)}
If $\al=0$ or [$\be_n=0$, for some $n\ge 0$], we are done;
so, without loss, one may assume that 
$\al> 0$ and $\be_n> 0$, $\forall n$.
As $\be_n\to 0< \al$, 
there must be some $n(\al)\ge 0$ in such a way that
$0< \al-\be_n< \al$, $\forall n\ge n(\al)$.
The imposed hypothesis now gives:
$\al-\be_n\in \Ga(y)$, $\forall n\ge n(\al)$.
Passing to limit as $n\to \oo$ yields
(by the semi-Archimedean property of $H$), $\al\in \Ga(y)$;
and the assertion follows.
\eproof

The following couple of properties will be useful in the sequel:

% Lemma 2
\blemma \label{le2}
The gauge  function $\ga$ is super-additive and subtractive:
\beq \label{305}
\ga(y_1+y_2)\ge \ga(y_1)+\ga(y_2),\ 
\mbox{if the right member exists}
\eeq
\beq \label{306}
\mbox{
$\ga(y_1-y_2)\le \ga(y_1)-\ga(y_2)$,\
whenever\ [$\ga(y_1)> -\oo$, $\ga(y_2)< \oo$].
}
\eeq
\elemma

\bproof
Without loss, one may assume that
$y_1,y_2\in H$ and
$\ga(y_1)> 0$, $\ga(y_2)> 0$.
By (\ref{302}), $y_1\ge k^0t_1, y_2\ge k^0t_2$,
whenever $0\le t_1< \ga(y_1), 0\le t_2< \ga(y_2)$;
and this yields (for all such $(t_1,t_2)$)
$y_1+y_2\ge k^0 [t_1+t_2]$\ (i.e.: $\ga(y_1+y_2)\ge t_1+t_2$).
This, and the arbitrariness of the precise couple, ends the argument.
The second part is directly obtainable from the first one,
in a standard way.
\eproof

In particular, when $Y$ is locally convex, (c02) holds provided
$H$ is closed. Then, our developments reduce to the ones in
Goepfert, Tammer and Z\u{a}linescu \cite{goepfert-tammer-zalinescu-2000}.
Note that, an axiomatic approach of these facts is possible,
under the lines in 
Artzner, Delbean, Eber and Heath  
\cite{artzner-delbean-eber-heath-1999};
we do not give details.

%  Section 4
\section{Main result}
\setcounter{equation}{0}

{\bf (A)}
Let in the following $Y$ stand for a (real) vector space.
Take a (non-degenerate, proper) Archimedean
convex cone $H$ of $Y$;
and let $(\le_H)$ stand for its induced quasi-order.
Further, let $K$ be some (non-degenerate, proper) semi-Archimedean
convex cone of $Y$, with $K\incl H$; 
and let $(\le_K)$ stand for the induced quasi-order. 

{\bf (B)}
Further, let $X$ be a nonempty set.
By a {\it pseudometric} over $X$  we mean any map
$d:X\times X\to R_+$; if, in addition, $d$ is
{\it reflexive} [$d(x,x)=0, \forall x\in X$]
and
{\it symmetric} [$d(x,y)=d(y,x), \forall x,y\in X$]
we say that it is a {\it rs-pseudometric}.
Let $(\Lb,\le)$ be some directed quasi-ordered structure.
Take a family $D=(d_\lb; \lb\in \Lb)$ of
rs-pseudometrics over $X$, with the properties: 
{\it $\Lb$-sufficient}
[$d_\lb(x,y)=0$, $\forall \lb\in \Lb$ $\limpl$ $x=y$],
$\Lb$-{\it monotone}
[$\lb\le \mu$ implies $d_\lb(.,.)\le d_\mu(.,.)$]
and $\Lb$-{\it triangular}
[$\forall \lb\in \Lb$, $\exists \mu\in \Lb(\lb,\le)$, with 
$d_\lb(x,z)\le d_\mu(x,y)+d_\mu(y,z)$, $\forall x,y,z\in X$].
By definition, $D$ will be referred to as a {\it Fang metric};
and $(X,D)$, as a {\it Fang uniform space}.

Technically speaking, $D$ may be viewed as 
a conv-Cauchy structure, in the following way.
Take an arbitrary sequence $(x_n; n\ge 0)$ in $X$.
Given $\lb\in \Lb$, 
the $d_\lb$-convergence of this sequence 
towards an $x\in X$ means: 
$d_\lb(x_n,x)\to 0$ as $n\to \oo$;
and it will be depicted as: $x_n\dlbtends x$.
If this holds for all $\lb\in \Lb$ then $(x_n; n\ge 0)$
is said to $D$-converge towards $x$; written as:
$x_n\Dtends x$; moreover, if $x\in X$ is generic in such 
a convention, $(x_n; n\ge 0)$ is called 
{\it $D$-convergent}.
On the other hand, given $\lb\in \Lb$, 
the $d_\lb$-Cauchy property of $(x_n; n\ge 0)$
means:
$\forall \veps> 0$, $\exists n:=n(\lb,\veps)$, such that 
$n\le p\le q$ $\limpl$ $d_\lb(x_p,x_q)< \veps$.
If this holds for each $\lb\in \Lb$, we say that $(x_n; n\ge 0)$ 
is $D$-{\it Cauchy}. 
Likewise, the $d_\lb$-semi-Cauchy property (where $\lb\in \Lb$)
of $(x_n; n\ge 0)$ is introduced as:
$\forall \veps> 0$, $\exists n:=n(\lb,\veps)$, such that 
$n\le m$  implies $d_\lb(x_n,x_m)< \veps$.
If this holds for all $\lb\in \Lb$, then we say that 
$(x_n; n\ge 0)$ is $D$-{\it semi-Cauchy}.
Note that, for each sequence in $X$,
\beq \label{401}
\barr{l}
\mbox{
$\forall \lb\in \Lb$:\
$d_\lb$-\mbox{Cauchy} $\limpl$ $d_\lb$-\mbox{semi-Cauchy;}  
}\\
\mbox{
hence:\ $D$-Cauchy $\limpl$ $D$-semi-Cauchy.
}
\earr
\eeq
The reciprocal of these is valid if {\it all} 
sequences in $X$ are involved; precisely, we have

% Lemma 3
\blemma \label{le3}
The global conditions below are equivalent to each other:
\bit
\item[(d01)]
$D$ is sequentially complete:\\
each $D$-Cauchy sequence is $D$-convergent
\item[(d02)]
$D$ is sequentially semi-complete:\\
each $D$-semi-Cauchy sequence is $D$-convergent.
\eit
\elemma

\bproof
By (\ref{401}), (d02) $\limpl$ (d01); 
so, it remains to prove that (d01) $\limpl$ (d02).
Let $(x_n; n\ge 0)$ be a $D$-semi-Cauchy sequence in $X$.
Fix $\lb\in \Lb$; and let $\mu\in \Lb(\lb,\le)$ 
given by the $\Lb$-triangular 
property of $D$. As $(x_n; n\ge 0)$ is 
$d_\mu$-semi-Cauchy, there exists, for the arbitrary fixed 
$\veps> 0$, some $n:=n(\mu,\veps)$ with:
$n\le m$ implies $d_\mu(x_n,x_m)< \veps/2$.
Combining these facts, yields
$$
n\le p\le q \limpl  
d_\lb(x_p,x_q)\le d_\mu(x_n,x_p)+d_\mu(x_n,x_q)< \veps;
$$
wherefrom, the desired assertion follows.
\eproof

{\bf (C)}
Now, let $(Y,H,K)$ be as above
and $(X,D)$ be a Fang uniform space.
Take a couple of functions function $F:X \to Y$, 
$k:\Lb\to K\sm (-H)$ with
\bit
\item[(d03)]
$F$ is $H$-bounded below:
$\exists b\in Y$ such that $G(x):=F(x)-b\in H$, $\forall x\in X$;
\item[(d04)]
$\lb \le \mu$ $\limpl$ $k(\lb)\le_K k(\mu)$\ \ 
($k(.)$ is increasing).
\eit
The relation $(\preceq_{(D,F)})$ over $X$ introduced as
\bit
\item[(d05)]
($x_1,x_2\in X$) $x_1\preceq_{(D,F)} x_2$ iff 
$k(\lb)d_\lb(x_1,x_2)\le_K F(x_1)-F(x_2)$,\ $\forall \lb\in \Lb$ 
\eit
is an order, as it can be directly seen.
For a number of both practical and theoretical reasons,
it would be useful to determine sufficient conditions under which 
$(\preceq_{(D,F)})$ is a Zorn order (cf. Section 2).
Essentially, these are 
\bit
\item[(d06)]
$D$ is sequentially $(K,F)$-complete:
each $D$-Cauchy sequence \\
$(x_n; n\ge 0)$ with $(F(x_n); n\ge 0)$, $K$-descending is 
$D$-convergent
\item[(d07)]
whenever $(x_n; n\ge 0)\incl X$ is $(\preceq_{(D,F)})$-ascending
and $x_n\Dtends x$ then $x_n \preceq_{(D,F)} x$, for each $n$. 
\eit
A basic particular case of (d07) is related to the property:
\bit
\item[(d08)]
$F$ is  sequentially $K$-descending $D$-lsc: \\
$x_n\Dtends x$ and $(F(x_n))$ is $K$-descending
imply $F(x_n)\ge_K F(x)$, $\forall n$.
\eit

Precisely, we have 

% Lemma 4
\blemma \label{le4}
Under the generally admitted facts (about $(H,K)$), (d08) $\limpl$ (d07).
\elemma 

\bproof
Let the sequence $(x_n; n\ge 0)$  in $X$ be 
$(\preceq_{(D,F)})$-ascending: 
\bit
\item[]
($\forall \lb\in \Lb$):\ 
$k(\lb)d_\lb(x_n,x_m)\le_K F(x_n)-F(x_m)$,\ \  
if $n\le m$;
\eit
clearly,  $(F(x_n))$ is $K$-descending.
Further, assume that $x_n\Dtends x$ for some $x\in X$.
By (d08), we have $F(x_n)\ge_K F(x)$, for all $n$;
this, by the working condition, yields
\beq \label{402}
(\forall \lb\in \Lb):\ 
k(\lb)d_\lb(x_n,x_m)\le_K F(x_n)-F(x),\ \  
\mbox{if}\ n\le m.
\eeq
Fix $\lb\in \Lb$ and $n\ge 0$. 
Let $\mu\in \Lb(\lb,\le)$
be the index assured by the
$\Lb$-triangular property of $D$.
From (\ref{402}) (and (d04))
\beq \label{403}
\barr{l}
k(\lb)d_\lb(x_n,x)\le_K 
k(\mu)d_\mu(x_n,x_m)+k(\lb)d_\mu(x_m,x) \\
\le_K F(x_n)-F(x)+ k(\lb)d_\mu(x_m,x),\
\mbox{for all}\ m\ge n;
\earr
\eeq
This, along with 
$k(\Lb)\incl K\sm (-K)$ and
the semi-Archimedean property of $K$, 
yields (via Lemma \ref{le1}) 
\beq \label{404}
(\forall \lb\in \Lb):\ 
k(\lb)d_\lb(x_n,x)\le_K F(x_n)-F(x);\ 
\mbox{i.e.}:\ x_n \preceq_{(D,F)} x.
\eeq
As $n\ge 0$ was arbitrarily fixed, (d07) holds.
\eproof

The main result of this exposition is

% Theorem 3
\btheorem \label{t3}
Let the (convex) cones $(H,K)$ in $Y$, 
the Fang metric $D$ and the 
couple of functions $(F,k)$ be such that 
(d03)-(d04) and (d06)-(d07) hold.
Then, for each $x_0\in X$ there exists $\bar x\in X$ with 

{\bf i)}
$x_0\preceq_{(D,F)} \bar x$;\ \ 
{\bf ii)}
$\bar x\preceq_{(D,F)} x'\in X$ $\limpl$  $\bar x=x'$.
\etheorem

Let us complete $Y$ with an element $\oo \notin Y$
and put $\wt Y=Y\cup \{\oo\}$;
the algebraic/order conventions involving this completion are
\bit
\item[]
($\oo=b+\oo=\oo+b$; $\oo=\lb \oo$;
$b\le_K \oo$, $\neg (\oo \le_K b)$),\
$\forall b\in Y$, $\forall \lb\in R_+^0$;
\eit
where $R_+^0:=]0,\oo[$.
As in Section 1, an "extended" form of this result 
is reached when $F:X\to \wt Y$ fulfills
\bit
\item[(d09)]
$F$ is proper:\  $\Dom(F):=\{x\in X; F(x)\ne \oo\}\ne \es$
\eit
as well (d03) (with $\Dom(F)$ in place of $X$).
But, this conclusion is obtainable from the above one: 
just apply Theorem \ref{t3} to 
the triplet $(X(x_0\preceq_{[D,F]});D;F)$, 
where $x_0\in \Dom(F)$ and $(\preceq_{[D,F]})$ 
is the quasi-order 
\bit
\item[(d10)]
($x_1,x_2\in X$):\ $x_1\preceq_{[D,F]} x_2$ iff 
$k(\lb) d_\lb(x_1,x_2)+F(x_2)\le_K F(x_1)$,\ $\forall \lb\in \Lb$; 
\eit
we do not give details.

\bproof {\bf (Theorem \ref{t3})}
There are several steps to be passed.
\sk

{\bf I)}
Fix $\theta\in \Lb$ 
and put $\Theta=\Lb(\theta,\le)$;
note that (as $(\Lb,\le)$ is directed),
\beq \label{405}
\mbox{
$\Theta$ is cofinal in $\Lb$:\ 
for each $\lb\in \Lb$ there exists $\mu\in \Theta$ with $\lb\le \mu$.
}
\eeq
Further, put $\de:=\ga_\theta$ 
[the gauge function attached to $(H;k(\theta))$] and 
$\psi=\de \circ G$.

{\bf II)}
Let $(\sqsubseteq_{(D,\psi)})$ stand for the relation:
\bit 
\item[(d11)]
($x_1,x_2\in X$):\ $x_1\sqsubseteq_{(D,\psi)} x_2$ iff 
$d_\lb(x_1,x_2)\le \psi(x_1)-\psi(x_2)$,\ $\forall \lb\in \Lb$; 
\eit
it is an order on $X$, as it can be directly seen. We show that 
\beq \label{406}
(\forall x_1,x_2\in X):\ 
x_1\preceq_{(D,F)} x_2 \limpl 
x_1\sqsubseteq_{(D,\psi)} x_2  
\limpl\ \psi(x_1)\ge \psi(x_2).
\eeq
The second part is clear; so, it remains to 
verify the first part.
Let $x_1,x_2\in X$ be such that 
$x_1\preceq_{(D,F)} x_2$; that is:
$$
k(\lb)d_\lb(x_1,x_2)\le F(x_1)-F(x_2)=G(x_1)-G(x_2),\ 
\forall \lb\in \Lb.
$$
As $k(.)$ is increasing, this yields
$$
k(\theta)d_\lb(x_1,x_2)\le G(x_1)-G(x_2),\ 
\forall \lb\in \Theta;
$$
so that, by Lemma \ref{le2} 
(and the imposed notations)
$$
d_\lb(x_1,x_2)\le \de (G(x_1)-G(x_2))\le
\psi(x_1)-\psi(x_2),\ 
\forall \lb\in \Theta.
$$
This, in turn, yields (as $\lb\mapsto d_\lb(.,.)$ is increasing
and (\ref{405}) holds)
$$
d_\lb(x_1,x_2)\le \psi(x_1)-\psi(x_2),\ 
\forall \lb\in \Lb; 
$$
that is: $x_1\sqsubseteq_{(D,\psi)} x_2$;  
hence the assertion.

{\bf III)}
We show that BB is applicable 
to $(X;\preceq_{(D,F)};\psi)$.
Firstly, by (\ref{406}), $\psi$ is 
decreasing (modulo $(\preceq_{(D,F)})$).
Secondly, let the sequence $(x_n; n\ge 0)$ in $X$ be
$(\preceq_{(D,F)})$-ascending:
\bit
\item[(d12)]
($\forall \lb\in \Lb$):\ 
$k(\lb)d_\lb(x_n,x_m)\le_K F(x_n)-F(x_m)$,\   
if $n\le m$.
\eit
note that, in such a case, 
$(F(x_n); n\ge 0)$ is $K$-descending.
By (\ref{406}), it follows that $(x_n; n\ge 0)$ 
is $(\sqsubseteq_{(D,\psi)}$)-ascending:
\beq \label{407}
(\forall \lb\in \Lb):\ 
d_\lb(x_n,x_m)\le \psi(x_n)-\psi(x_m),\ 
\mbox{if}\ n\le m;
\eeq
and, from this, $(x_n; n\ge 0)$ is $D$-Cauchy.
Combining with  (d06),  it follows that
$x_n\Dtends x$, for some $x\in X$. Moreover,
let $n\ge 0$ be arbitrary fixed. By (d07), 
\beq \label{408}
(\forall \lb\in \Lb):\ 
k(\lb)d_\lb(x_n,x)\le_K F(x_n)-F(x);\ \  
\mbox{i.e.:}\ x_n\preceq_{(D,F)} x.
\eeq
This shows that 
$(X,\preceq_{(D,F)})$ is sequentially inductive;
and proves the claim.

Applying BB to these data, one gets 
that, for $x_0\in X$, there exists $\bar x\in X$ with

{\bf j)}
$x_0 \preceq_{(D,F)} \bar x$;\ \ 
{\bf jj)}
$\bar x \preceq_{(D,F)} x'\in X$ 
$\limpl$  $\psi(\bar x)=\psi(x')$.

\n
The former of these is just {\bf i)}.
Moreover, by the latter of these, one gets {\bf ii)}.
For, let $x'\in X$ be such that 
$\bar x\preceq_{(D,F)} x'$. 
Again by (\ref{406}), $\bar x\preceq_{(D,\psi)} x'$; 
and, by {\bf jj)} above,  $\psi(\bar x)=\psi(x')$.
Combining these gives $\bar x=x'$;
hence the assertion. 
\eproof

In particular, when $Y$ is a locally convex space,
the Archimedean property of $H$ 
is assured when $H=\cl(K)$; 
moreover, (d06) holds under (d01)/(d02).
The corresponding version of Theorem \ref{t3}
is just the main result in 
Zhu and Li \cite{zhu-li-2007}
proved via rather different methods.
On the other hand, 
when $k(.)$ is a constant function, 
we get the main result in
Turinici \cite{turinici-2002};
which includes the ones in 
Goepfert, Tammer and Z\u{a}linescu 
\cite{goepfert-tammer-zalinescu-2000}.
But, as precise by the quoted authors, 
their statements include EVP; 
hence, so does Theorem \ref{t3}.
Summing up, we must have
(DC $\limpl$) BB $\limpl$ Theorem \ref{t3} $\limpl$ EVP;
wherefrom, by the developments of Section 2, 
Theorem \ref{t3} is equivalent with both BB and EVP.

%  Section 5
\section{Scalar versions}
\setcounter{equation}{0}

Let $Y$ be a (real) vector space.
By the developments in Lemma \ref{le4},
the choice 
\bit
\item[]
$H=K$=Archimedean (non-degenerate, proper) 
(convex) cone of $Y$
\eit
is allowed in Theorem \ref{t3}.
This, in the case of $Y=R$, $H=K=R_+$, 
yields a  variational principle 
over Fang uniform spaces, including 
Hamel's \cite{hamel-2005}.
It is our aim in the following to 
state this principle; 
as well as to discus a lot of related facts.

Let $X$ be a nonempty set; and
$(\Lb,\le)$ be a directed quasi-ordered structure.
Take a family $D=(d_\lb; \lb\in \Lb)$ of
rs-pseudometrics over $X$, with the properties: 
{\it $\Lb$-sufficient},
$\Lb$-{\it monotone}
and $\Lb$-{\it triangular};
by a previous convention, $D$ will be referred to as a 
Fang metric. 
Define a conv-Cauchy structure on $(X,D)$ 
as in Section 4. 
Further, let ($\vphi:X\to R$; $h:\Lb\to R_+^0$) 
be a couple of functions with the properties (a01) and
\bit
\item[(e01)]
$\lb\le \mu$ $\limpl$ $h(\lb)\le h(\mu)$\ \ 
($h(.)$ is increasing).
\eit
The relation $(\preceq_{(D,\vphi)})$ over $X$ introduced as
\bit
\item[(e02)]
($x_1,x_2\in X$):\ $x_1\preceq_{(D,\vphi)} x_2$ iff 
$h(\lb)d_\lb(x_1,x_2)\le \vphi(x_1)-\vphi(x_2)$,\ $\forall \lb\in \Lb$ 
\eit
is an order, as it can be directly seen;
as in Section 4, we want to determine
sufficient conditions under which 
$(\preceq_{(D,\vphi)})$ be a Zorn one.
As precise there (cf. Lemma \ref{le4}), 
the specific assumptions to be added write
\bit
\item[(e03)]
$D$ is sequentially $\vphi$-complete:
each $D$-Cauchy sequence \\
$(x_n; n\ge 0)$  with $(\vphi(x_n); n\ge 0)$, descending 
is $D$-convergent
\item[(e04)]
$\vphi$ is  sequentially descending $D$-lsc: \\
$\lim_n \vphi(x_n)\ge \vphi(x)$, whenever $x_n\Dtends x$
and $(\vphi(x_n))$ is descending.
\eit

The appropriate answer to this question is contained in

% Theorem 4
\btheorem  \label{t4}
Let the Fang metric $D$ and the 
functions $[\vphi; h]$ be as in (a01), (e01) and (e03)-(e04).
Then, for each $u\in X$, there exists $v\in X$ with
\beq \label{501}
\mbox{
$h(\lb)d_\lb(u,v)\le \vphi(u)-\vphi(v)$, $\forall \lb\in \Lb$\
(hence $\vphi(u)\ge \vphi(v)$)
}
\eeq
\beq \label{502}
\forall x\in X\sm \{v\}, \exists \mu=\mu(x)\in \Lb:
h(\mu)d_\mu(v,x)> \vphi(v)-\vphi(x).
\eeq
\etheorem

In particular, (e03) holds under (d01)/(d02);
when Theorem \ref{t4} is just 
Hamel's variational principle
\cite{hamel-2005} 
(in short: HVP).
This last result -- 
based on a maximal principle comparable with
Br{\o}ndsted's \cite{brondsted-1974} --
extends the related statement in 
Fang \cite{fang-1996},
obtained {\it via} Zorn maximal techniques.
It also includes the contribution due to
Had\v{z}i\'{c} and \v{Z}iki\'{c} \cite{hadzic-zikic-1998}
(in short: HZVP), founded on the maximal principle in 
Hicks \cite{hicks-1989}; 
we do not give details.
\sk

Now, Theorem \ref{t4} is but a particular version of Theorem \ref{t3};
hence, it is reducible to BB.
On the other hand, HVP includes EVP;
just take $\Lb$ as a singleton.
Summing up, 
BB $\limpl$ Theorem \ref{t4} $\limpl$ HVP $\limpl$ EVP;
this, by the developments of Section 2, tells us that 
HVP is equivalent with both BB and EVP.

Concerning the former of these inclusions
(BB $\limpl$ Theorem \ref{t4}) it is worth noting that,
in Theorem \ref{t3}, BB was applied in a "local" way,
by means of the point $\theta\in \Lb$ and its attached
section $\Theta:=\Lb(\theta,\le)$.
However, the presence of a "scalar" objective function $\vphi$
(in place of the vectorial function $F$)
suggests us that a "global" application of BB is highly expectable.
To see that this is indeed the case, 
it would be useful working with the 
particular choice of $h(.)$, taken as
\bit
\item[(e05)]
$h(\lb)=1$, for all $\lb \in \Lb$ (so, $h(.)$ is constant).
\eit
This yields the following "standard" version of
Theorem \ref{t4} to be considered:

% Theorem 5
\btheorem  \label{t5}
Let the Fang metric $D$ and the 
function $\vphi$ be as in (a01) and (e03)-(e04).
Then, for each $u\in X$, there exists $v\in X$ with
\beq \label{503}
\mbox{
$d_\lb(u,v)\le \vphi(u)-\vphi(v)$, $\forall \lb\in \Lb$\
(hence $\vphi(u)\ge \vphi(v)$)
}
\eeq
\beq \label{504}
\forall x\in X\sm \{v\}, \exists \mu=\mu(x)\in \Lb:
d_\mu(v,x)> \vphi(v)-\vphi(x).
\eeq
\etheorem

Concerning the relationships between these results,
the following answer holds:

% Proposition 7
\bprop \label{p7}
Under these general assumptions, we have
Theorem \ref{t5} $\limpl$ Theorem \ref{t4};
hence Theorem \ref{t5} $\lequi$ Theorem \ref{t4}.
\eprop

\bproof
Let the conditions of Theorem \ref{t4} be admitted.
Define another family $E=(e_\lb; \lb\in \Lb)$ of
rs-pseudometrics over $X$ according to
\bit
\item[(e06)]
$e_\lb(x,y)=h(\lb)d_\lb(x,y)$, $x,y\in X$.
\eit
The $\Lb$-sufficiency of $E$ results at once from that
of $D$; and the $\Lb$-monotonicity of the same
is directly reducible to that of $h(.)$.
Finally, we claim that $E$ is $\Lb$-triangular.
Let $\lb\in \Lb$ be arbitrarily fixed.
By the $\Lb$-triangular property of $D$,
there exists $\mu\in \Lb(\lb,\le)$ such that 
$d_\lb(x,z)\le d_\mu(x,y)+d_\mu(y,z)$, $\forall x,y,z\in X$.
This, again by the increasing property of
$h(.)$, yields
$$ \barr{l}
e_\lb(x,z)=h(\lb)d_\lb(x,z)\le
h(\lb)[d_\mu(x,y)+d_\mu(y,z)] \\
\le e_\mu(x,y)+e_\mu(y,z),\ \forall x,y,z\in X;
\earr
$$
and the assertion follows.
Summing up, $E$ is a Fang metric;
it may generate a conv-Cauchy structure  on $X$, 
by the construction in Section 4.
Concerning its connections with the 
Fang metric $D$ (and its attached conv-Cauchy structure), 
one has (for all sequences $(x_n)$ in $X$, and all $x\in X$)
\beq \label{505}
[\forall \lb\in \Lb:
(x_n \dlbtends x) \lequi (x_n \elbtends x)];\
\mbox{hence}\ [(x_n \Dtends x) \lequi  (x_n \Etends x)];
\eeq
as well as (for a generic sequence $(x_n)$ in $X$)
\beq \label{506}
\mbox{
[$\forall \lb\in \Lb$: 
$d_\lb$-Cauchy $\lequi$ $e_\lb$-Cauchy];\
hence $D$-Cauchy $\lequi$ $E$-Cauchy.
}
\eeq
The conv-Cauchy structures attached to 
the Fang metrics $D$ and $E$ 
are thus equivalent to each other.
As a direct consequence,
(e03)-(e04) are holding over $(E,\vphi)$; and we are done.
\eproof

Having this precise, we may now return to 
the addressed question. 

\bproof {\bf (Theorem \ref{t5})}
Let the conditions of Theorem \ref{t5} hold;
and $(\preceq_{(D,\vphi)})$ stand for the order (e02),
where $h(.)$ is taken as in (e05).
We have to verify that BB is applicable to 
$(X,\preceq_{(D,\vphi)};\vphi)$.
Clearly, $\vphi$ is descending [modulo $(\preceq_{(D,\vphi)})$].
Moreover, let $(x_n; n\ge 0))$ be an ascending 
[modulo $(\preceq_{(D,\vphi)})$] sequence in $X$:
\bit
\item[(e07)]
($\forall \lb\in \Lb$):\ 
$d_\lb(x_n,x_m)\le \vphi(x_n)-\vphi(x_m)$,\ \  
if $n\le m$.
\eit
The sequence $(\vphi(x_n))$ is descending and bounded
from below; hence a Cauchy one
[$\forall \de> 0, \exists n(\de):\
n(\de)\le p\le q\limpl \vphi(x_p)-\vphi(x_q)< \de$].
This, along with (e07), tells us that $(x_n)$ is $D$-Cauchy.
Taking (e03) into account, it follows that $x_n \Dtends x$ 
as $n\to \oo$, for some $x\in X$.
Combining with (e04) and Lemma \ref{le4} gives 
$x_n\preceq_{(D,\vphi)} x$, $\forall n$; 
wherefrom, $(X,\preceq_{(D,\vphi)})$ is sequentially inductive.
From BB it then follows that, for the starting $u\in X$, 
there exists a $(\preceq_{(D,\vphi)},\vphi)$-maximal $v\in X$ with 
$u\preceq_{(D,\vphi)} v$.
This element has the properties (\ref{503})+(\ref{504}), 
and the conclusion follows. 
\eproof

An alternate argument for establishing this result 
is the following.
Let the Fang uniform space $(X,D)$ 
and the function $\vphi:X\to R$ be as in Theorem \ref{t5}.
The associated family of relations
$\calv=\{U(\lb,r); \lb\in \Lb, r> 0\}$ given by 
\bit
\item[(e08)]
$U(\lb,r)=\{(x,y)\in X\times X; e_\lb(x,y)< r\}$,
$\lb\in \Lb, r> 0$,
\eit
is a fundamental system of entourages for a uniform structure
$\calu=\calu(\Lb,D)$ over $X$
in the  sense described by
Bourbaki \cite[Ch 2, Sect 1]{bourbaki-1989}).
Further, let the ordering $(\preceq_{(D,\vphi)})$ over $X$ 
be defined as in (e02).
We claim that BMP (Proposition \ref{p5}) is applicable to these data.
This will follow from the

\bproof {\bf [(alternate) Theorem \ref{t5}]}
Clearly, $\calv$ is (sequentially) $(\preceq_{(D,\vphi)})$-complete,
if we take (e03) into account. Moreover, 
$\calv$ is $(\le,\vphi)$-admissible, {\it via} 
\beq \label{507}
(\forall \de> 0):\  [x\le y,\ \vphi(x)-\vphi(y)< \de] \limpl 
(x,y)\in \cap\{U(\lb,\de); \lb\in \Lb\}.
\eeq
It remains to verify that $(\preceq_{D,\vphi)})$ is $\calv$-selfclosed. 
Let $(x_n; n\ge 0)$ in $X$ be ascending (cf. (e07));
with, in addition, $x_n\Dtends x$.
Note that, in such a case, $(\vphi(x_n))$ is descending;
so that $\vphi(x_n)\ge \vphi(x)$, $\forall n$,
if one takes (e04) into account.
Given $\lb\in \Lb$, take $\mu\in \Lb(\lb,\le)$ according to the 
$\Lb$-triangular property of $D$; 
and let $n$ be arbitrary fixed.
By (e04) and the remark above
$$
d_\lb(x_n,x)\le d_\mu(x_n,x_m)+d_\mu(x_m,x)\le 
\vphi(x_n)-\vphi(x)+d_\mu(x_m,x),\ \forall m\ge n. 
$$
Passing to limit as $m\to \oo$ we get
$d_\lb(x_n,x)\le \vphi(x_n)-\vphi(x)$, $\forall \lb\in \Lb$
[that is: $x_n\preceq_{(D,\vphi)} x$], for all $n$; hence the claim.
\eproof

The following completion of these facts is to be noted.
Let $I$ be some nonempty set.
Take a family $F=(f_i; i\in I)$ of
rs-pseudometrics over $X$, with the properties:  
{\it $I$-sufficient}
[$f_i(x,y)=0$, for all  $i\in I$ imply $x=y$],
and {\it $I$-triangular}
[for each $i\in I$, there exist $j=j(i)$  and $k=k(i)$ in $I$ 
such that
$f_i(x,z)\le f_j(x,y)+f_k(y,z)$, $\forall x,y,z\in X$].
In this case, the couple $(X,F)$ will be termed a
{\it BMLO uniform space};
see Benbrik, Mbarki, Lahrech and Ouahab 
\cite{benbrik-mbarki-lahrech-ouahab-2006}.
Clearly, any Fang uniform space is a BMLO uniform space as well.
But, the reciprocal inclusion is also true.
In fact, let $\Lb$ stand for the class of all
(nonempty) finite parts of $I$, endowed with 
the usual inclusion, $(\incl)$;
note that, $(\Lb,\incl)$ is a directed ordered structure.
For each $\lb\in \Lb$ define the rs-pseudometric
$d_\lb$ over $X$ as:
$d_\lb(x,y)=\sup\{f_i(x,y); i\in \lb\}$, $x,y\in X$.
The family $D=(d_\lb; \lb\in \Lb)$ of all these 
is easily shown to be $\Lb$-sufficient, 
$\Lb$-monotone and $\Lb$-triangular; 
i.e., $(X;D)$ is a Fang uniform space.
In addition, all usual $F$-concepts 
(like $F$-convergence and $F$-Cauchy) are 
equivalent to their corresponding $D$-concepts.
Hence, all variational results over BMLO uniform spaces
established by these authors are completely reducible to 
those involving Fang uniform spaces we just presented; see also
Hamel and Loehne \cite{hamel-loehne-2003}.
In particular, this is retainable for the 
variational principles in standard uniform spaces due to 
Mizoguchi \cite{mizoguchi-1990},
because any such structure is a BMLO uniform space.
Further aspects may be found in
Had\v{z}i\'{c} and Ovcin \cite{hadzic-ovcin-1994};
see also
Chang et al \cite{chang-cho-lee-jung-kang-1997}.

%  Section 6
\section{(ZF-AC) approach}
\setcounter{equation}{0}

Let us now return to the setting of Theorem \ref{t5}.
Precisely, let $X$ be a nonempty set; 
and $(\Lb,\le)$, a directed quasi-ordered structure.
Further, take a family $D=(d_\lb; \lb\in \Lb)$
of rs-pseudometrics over $X$, 
with the properties:
$\Lb$-sufficient, 
$\Lb$-monotone 
and 
$\Lb$-triangular.
We term it, a {\it Fang metric};
note that, by the developments in Section 4, 
$D$ may generate a conv-Cauchy structure on $X$. 
Further, let $\vphi:X\to R$ be a function
as in (a01); and let the couple $(D,\vphi)$ fulfill
(e03)-(e04). Then (cf. Section 5), 
conclusions of Theorem \ref{t5} are holding for our data.

In particular, when $\Lb$ (hence $D$ as well) is a
singleton, Theorem \ref{t5} yields the following
metric variational statement.
Let $d$ be a metric on $X$; and $\vphi:X\to R$
be some function as in (a01).

% Theorem 6
\btheorem \label{t6}
Suppose that, in addition,
\bit
\item[(f01)]
$d$ is $\vphi$-complete: \\
each $d$-Cauchy sequence 
with $(\vphi(x_n))$ descending is $d$-convergent
\item[(f02)]
$\vphi$ is descending $d$-lsc: \\
$\lim_n \vphi(x_n)\ge \vphi(x)$, whenever $x_n\dtends x$
and $(\vphi(x_n))$ is descending.
\eit
Then, conclusions of EVP are holding.
\etheorem

Combining with the previous facts, we have:
BB $\limpl$ Theorem \ref{t5} $\limpl$ Theorem \ref{t6} $\limpl$ EVP.
This, by the developments in Section 2, yields:
Theorem \ref{t5} and Theorem \ref{t6} are equivalent with both 
BB and EVP; hence, mutually equivalent.

Generally, the family of rs-pseudometrics 
$D=(d_\lb; \lb\in \Lb)$ is non-denumerable.
For example, in case of the Fang uniform spaces
constructed from a probabilistic metric space
(cf. Fang \cite{fang-1996})
or fuzzy metric spaces
(cf. Had\v{z}i\'{c} and \v{Z}iki\'{c} \cite{hadzic-zikic-1998}) 
we have $(\Lb,\le):=(]0,1],\ge)$; here $(\ge)$ is the usual 
dual ordering on $R$.
On the other hand, the particular (modulo (e05))
ordering (e02) appearing there
\bit
\item[(f03)]
($x_1,x_2\in X$):\ $x_1\preceq_{(D,\vphi)} x_2$ iff 
$d_\lb(x_1,x_2)\le \vphi(x_1)-\vphi(x_2)$,\ $\forall \lb\in \Lb$ 
\eit
may be ultimately viewed as a Br{\o}ndsted one, by simply taking the 
supremum in the left hand side of this relation.
So, we may ask whether a deduction of 
Theorem \ref{t5} from Theorem \ref{t6} is possible.
The (positive) answer to this is contained in

% Proposition 8
\bprop \label{p8}
We have, in (ZF-AC) (without any use of (DC)):
\beq \label{601}
\mbox{
Theorem \ref{t6} $\limpl$ Theorem \ref{t5}
[hence Theorem \ref{t6} $\lequi$ Theorem \ref{t5}].
}
\eeq
\eprop

Before passing to the effective part, we need
some preliminary facts.
Denote
\bit
\item[(f04)]
$\De(x,y)=\sup\{d_\lb(x,y); \lb\in \Lb\}$, $x,y\in X$.
\eit
Since all members of $D$ are rs-pseudometrics,
$\De$ is also endowed with such properties.
Moreover (as $D$ is $\Lb$-triangular),
$\De$ is
{\it triangular}
[$\De(x,z)\le \De(x,y)+\De(y,z), \forall x,y,z\in X$];
finally, $\De$ is
{\it sufficient} [$\De(x,y)=0 \limpl x=y$]; because so is $D$.
Summing up, $\De$ is a generalized metric on $X$,
in the Luxemburg-Jung sense
\cite{luxemburg-1958}, \cite{jung-1969}.
It allows us introducing a 
conv-Cauchy structure on $X$ as
\bit
\item[(f05)]
$x_n\Detends x$ iff $\De(x_n,x)\to 0$ as $n\to \oo$
\item[(f06)]
$(x_n)$ is $\De$-Cauchy iff $\lim_{n,m}\De(x_n,x_m)=0$.
\eit
Alternatively, it allows us introducing 
a uniform  structure $\calu=\calu(\De)$ on $X$
as the one for which
$\calv=\{U(\veps); \veps> 0\}$, where
\bit
\item[(f07)]
$U(\veps)=\{(x,y\in M\times M; \De(x,y)< \veps\}$, $\veps> 0$,
\eit
is a fundamental system of entourages.
The natural question to be posed is that of clarifying
the relationships between these and those 
attached to the family $D=(d_\lb; \lb\in \Lb)$.
First, by these conventions, we have

% Lemma 5
\blemma \label{le5}
The generic local inclusions hold:
\beq \label{602}
(\forall (x_n), \forall x):\ [x_n \Detends x] \limpl [x_n \Dtends x].
\eeq
\beq \label{603}
\mbox{
(for each sequence):\ $\De$-Cauchy $\limpl$ $D$-Cauchy.
}
\eeq
\elemma

The reciprocal inclusions are not in general true; 
because the uniform structure attached to $D$ is strictly finer
than that induced by the generalized metric $\De$.
A useful completion of these facts is contained in

% Lemma 6
\blemma \label{le6}
Under these notations,
\beq \label{604}
\mbox{
($\forall (x_n),\forall x$)
[$(x_n)$ is $\De$-Cauchy,\ $x_n \Dtends x]$ imply $[x_n \Detends x]$
}
\eeq
\beq \label{605}
\mbox{
$D$ is sequentially $\vphi$-complete $\limpl$ $\De$ is $\vphi$-complete.
}
\eeq
\elemma

\bproof ({\bf Lemma \ref{le6}})
The second part in the statement follows ({\it via} Lemma \ref{le5}) 
from the first part of the same; 
so, it is sufficient proving that
(\ref{604}) holds.
Let $(x_n)$ be a $\De$-Cauchy sequence in $X$, so as (for some $x\in X$)
$$  \mbox{
$x_n \Dtends x$\ (hence $d_\lb(x_n,x) \to 0$, for each $\lb\in \Lb$).
}
$$
By definition, for each $\be> 0$ there exists some rank $n(\be)$
in such a way that
$\De(x_i,x_j)\le \be$ (hence $d_\lb(x_i,x_j)\le \be$, $\forall \lb\in \Lb$),
whenever $n(\be)\le i\le j$.
Let the rank $i\ge n(\be)$ be arbitrarily fixed,
and, for each $\lb\in \Lb$, let $\mu\in \Lb(\lb,\le)$
be the index given by the $\Lb$-triangular property of $D$.
We have, for all such $(\lb,\mu)$,
$$
d_\lb(x_i,x)\le d_\mu(x_i,x_j)+d_\mu(x_j,x)\le
\be +d_\mu(x_j,x),\  \forall j\ge i.
$$
Passing to limit upon $j$ gives (for all $i$ like before)
$$ \mbox{
$d_\lb(x_i,x)\le \be,\ \forall \lb\in \Lb$\
(hence $\De(x_i,x)\le \be$).
}
$$
This, by the arbitrariness of $\be$, yields $x_n \Detends x$;
as claimed.
\eproof

Having these precise, we may now pass to the effective argument.

\bproof {\bf (Proposition \ref{p8})}
Let $\De$ stand for the generalized metric over $X$
introduced {\it via} (f04); and put
$X_u=\{x\in X; \De(u,x)\le \vphi(u)-\vphi(x)\}$
(where $u\in X$ is the point in Theorem \ref{t5}).
Clearly, $\De$ is a standard metric over $X_u$.
On the other hand, the imposed conditions assure us
({\it via} Lemma \ref{le5} and Lemma \ref{le6}) 
that (f01)+(f02) hold (over $X_u$) with $\De$ instead of $d$.
Summing up,
Theorem \ref{t6} applies to $(X_u;\De;\vphi)$.
It gives us, for the starting $u\in  X_u$
some $v\in X_u$
fulfilling (\ref{101})+(\ref{102}) (relative to $X_u$ and $\De$).
This yields the conclusions (\ref{503})+(\ref{504}) we want.
\eproof

Summing up, (DC) is not needed to establish the 
logical equivalence in Proposition \ref{p8}.
However, both Theorem \ref{t5} and Theorem \ref{t6}
are deductible in (ZF-AC+DC). 
Further aspects will be delineated elsewhere.

%  References


\begin{thebibliography}{99}


% 1
\bibitem{artzner-delbean-eber-heath-1999}
{P. Artzner}, {F. Delbean}, {J. M. Eber} and  {D. Heath}, 
\it Coherent measures of risk,
\rm Math. Finance, 9 (1999), 203-228.


% 2
\bibitem{benbrik-mbarki-lahrech-ouahab-2006}
{A. Benbrik}, {A. Mbarki}, {S. Lahrech} and {A. Ouahab},
\it Ekeland's principle for vector-valued maps based on the characterization 
of uniform spaces via families of generalized quasi-metrics,
\rm Lobachevskii J. Math., 21 (2006), 33-44.


% 3
\bibitem{bernays-1942}
{P. Bernays},
\it A system of axiomatic set theory: Part III. Infinity and enumerability analysis,
\rm J. Symbolic Logic, 7 (1942), 65-89.


% 4
\bibitem{bourbaki-1989}
{N. Bourbaki},
\it General Topology (Chs 5-10),
\rm Springer, Berlin, 1989.


% 5
\bibitem{brezis-browder-1976}
{H. Brezis} and {F. E. Browder},
\it A general principle on ordered sets in nonlinear functional analysis,
\rm Advances Math., 21 (1976), 355-364.


% 6
\bibitem{brondsted-1974}
{A. Br{\o}ndsted},
\it  On a lemma of Bishop and Phelps,
\rm  Pacific J. Math., 55 (1974), 335-341.


% 7
\bibitem{brondsted-1976}
{A. Br{\o}ndsted},
\it Fixed points and partial orders,
\rm Proc. Amer. Math. Soc., 60 (1976), 365-366.


% 8
\bibitem{brunner-1987}
{N. Brunner},
\it Topologische Maximalprinzipien,
\rm Zeitschr. Math. Logik Grundl. Math., 33 (1987), 135-139.


% 9
\bibitem{chang-cho-lee-jung-kang-1997}
{S. S. Chang}, {Y. J. Cho}, {B. S. Lee}, {J. S. Jung} and {S. M. Kang},
\it Coincidence point theorems and minimization theorems in fuzzy metric spaces,
\rm  Fuzzy Sets Syst., 88 (1997), 119-127.


% 10
\bibitem{cristescu-1977}
{R. Cristescu},
\it Topological Vector Spaces,
\rm  Noordhoff Intl. Publishers, Leyden (The Netherlands), 1977.


% 11
\bibitem{ekeland-1979}
{I. Ekeland},
\it Nonconvex minimization problems,
\rm Bull. Amer. Math. Soc. (New Series), 1 (1979), 443-474.


% 12
\bibitem{fang-1996}
{J. X. Fang},
\it The variational principle and fixed point theorems in certain topological spaces,
\rm J. Math. Anal. Appl., 202 (1996), 398-412.


% 13
\bibitem{goepfert-tammer-zalinescu-2000}
{A. Goepfert}, {C. Tammer} and {C. Z\u{a}linescu},
\it On the vectorial Ekeland's variational principle and minimal points in product spaces,
\rm Nonlin. Anal., 39 (2000), 909-922.


% 14
\bibitem{hadzic-ovcin-1994}
{O. Had\v{z}i\'{c}} and {Z. Ovcin},
\it Fixed point theorem in fuzzy metric spaces and probabilistic metric spaces,
\rm Review Res. Fac. Sci. Novi Sad Univ. (Math. Series), 24 (1994), 197-209.


% 15
\bibitem{hadzic-zikic-1998}
{O. Had\v{z}i\'{c}} and {T. \v{Z}iki\'{c}},
\it On Caristi's fixed point theorem in F-type topological spaces,
\rm Novi Sad J. Math., 28 (1998), 91-98.


% 16
\bibitem{hamel-2005}
{A. Hamel},
\it Equivalents to Ekeland's variational principle in uniform spaces,
\rm Nonlin. Anal., 62 (2005), 913-924.


% 17
\bibitem{hamel-loehne-2003}
{A. Hamel} and {A. Loehne},
\it A minimal point theorem in uniform spaces,
\rm in "Nonlinear Analysis and Applications: To V. Lakshmikantham on his 80th birthday"
[R. P. Agarwal and D. O'Regan, eds.], vol. 1, pp. 577-593, Kluwer, Dordrecht, 2003.


% 18
\bibitem{hicks-1989}
{T. L. Hicks},
\it Some fixed point theorems,
\rm Radovi Mat., 5 (1989), 115-119.


% 19
\bibitem{hyers-isac-rassias-97}
{D. H. Hyers}, {G. Isac} and {T. M. Rassias},
\it Topics in Nonlinear Analysis and Applications,
\rm World Sci. Publ., Singapore, 1997.


% 20
\bibitem{jung-1969}
{C. F. K. Jung},
\it On generalized complete metric spaces,
\rm Bull. Amer. Math. Soc., 75 (1969), 113-116.


% 21
\bibitem{kang-park-1990}
{B. G. Kang} and {S. Park},
\it On generalized ordering principles in nonlinear analysis,
\rm Nonlin. Anal., 14 (1990), 159-165.


% 22
\bibitem{kasahara-1976}
{S. Kasahara},
\it On some generalizations of the Banach contraction theorem,
\rm Publ. Res. Inst. Math. Sci. Kyoto Univ., 12 (1976), 427-437.


% 23
\bibitem{luxemburg-1958}
{W. A. J. Luxemburg},
\it On the convergence of successive approximations 
in the theory of ordinary differential equations (II),
\rm Indagationes Math., 20 (1958), 540-546.


% 24
\bibitem{mizoguchi-1990}
{N. Mizoguchi},
\it A generalization of Br{\o}ndsted's result and its applications,
\rm Proc. Amer. Math. Soc., 108 (1990), 707-714.


% 25
\bibitem{moore-1982}
{G. H. Moore},
\it Zermelo's Axiom of Choice: its Origin, Development and Influence,
\rm Springer, New York, 1982.


% 26
\bibitem{nachbin-1965}
{L. Nachbin},
\it Topology and Order,
\rm D. van Nostrand Comp. Inc., Princeton (N.J.), 1965.


% 27
\bibitem{schechter-1997}
{E. Schechter},
\it Handbook of Analysis and its Foundation,
\rm Academic Press, New York, 1997. 


% 28
\bibitem{tarski-1948}
{A. Tarski},
\it Axiomatic and algebraic aspects of two theorems on sums of cardinals,
\rm Fund. Math., 35 (1948), 79-104.


% 29
\bibitem{turinici-1994}
{M. Turinici},
\it Vector extensions of the variational Ekeland's result,
\rm An. \c{S}t. Univ. "A. I. Cuza" Ia\c{s}i (S I-a: Mat), 40 (1994), 225-266.


% 30
\bibitem{turinici-2002}
{M. Turinici},
\it Minimal points in product spaces,
\rm An. \c{S}t. Univ. "Ovidius" Constan\c{t}a (Ser. Math.), 10 (2002), 109-122.


% 31
\bibitem{turinici-2011}
{M. Turinici},
\it  Brezis-Browder principle and Dependent Choice,
\rm An. \c{S}t. Univ. "Al. I. Cuza" Ia\c{s}i (S. N.), Mat., 57 (2011), 263-277.


% 32
\bibitem{wolk-1983}
{E. S. Wolk},
\it On the principle of dependent choices and some forms of Zorn's lemma,
\rm Canad. Math. Bull., 26 (1983), 365-367.


% 33
\bibitem{zhu-li-2007}
{J. Zhu} and {S. J. Li},
\it Generalization of ordering principles and applications,
\rm J. Optim. Th. Appl., 132 (2007), 493-507.


\end{thebibliography}
\end{document}